\documentclass[10pt]{amsart}

\usepackage{amsmath}
\usepackage{amssymb}
\usepackage{mathtools}
\usepackage{fancyhdr}
\usepackage{stmaryrd}
\usepackage[shortlabels]{enumitem}
\usepackage{hyperref}
\hypersetup{
	hidelinks,
	pdftitle={The D-property in domains and monoids},
	pdfauthor={Mingyang Cen, Felix Gotti, Eric Guo, Patrek Ragnarsson}
}

\newtheorem{theorem}{Theorem}[section]
\newtheorem{prop}[theorem]{Proposition}
\newtheorem{lemma}[theorem]{Lemma}
\newtheorem{corollary}[theorem]{Corollary}
\theoremstyle{definition}
\newtheorem{definition}[theorem]{Definition}
\newtheorem{remark}[theorem]{Remark}
\newtheorem{example}[theorem]{Example}
\newtheorem{question}[theorem]{Question}
\numberwithin{equation}{section}

\newcommand{\nn}{\mathbb{N}}
\newcommand{\qq}{\mathbb{Q}}
\newcommand{\rr}{\mathbb{R}}
\newcommand{\zz}{\mathbb{Z}}
\newcommand{\gp}{\mathcal{G}}
\providecommand\ldb{\llbracket}
\providecommand\rdb{\rrbracket}
\DeclareMathOperator{\Frac}{Frac}
\DeclareMathOperator{\supp}{supp}
\newcommand{\PSP}{\mathrm{PSP}}
\setlist[enumerate]{label=\textup{(\arabic*)},leftmargin=*}
\setlist[itemize]{leftmargin=*}
\title{The D-property in domains and monoids}

\author{Mingyang Cen}
\address{PRIMES-USA\\MIT\\Cambridge, MA 02139}
\email{mingyangcen@gmail.com}

\author{Felix Gotti}
\address{Department of Mathematics\\MIT\\Cambridge, MA 02139}
\email{fgotti@mit.edu}

\author{Eric Guo}
\address{PRIMES-USA\\MIT\\Cambridge, MA 02139}
\email{ericguo01@gmail.com}

\author{Patrek Ragnarsson}
\address{PRIMES-USA\\MIT\\Cambridge, MA 02139}
\email{patrek.karason@gmail.com}

\date{November 25, 2026}


\begin{document}
	
	\begin{abstract}
		An integral domain satisfies the D-property if $a \mid bc$ and coprimality of $a$ and $b$ imply $a \mid c$ for nonzero elements $a$, $b$, and $c$. For a domain $R$ with fraction field $K$ and a nontrivial Puiseux monoid $M$, we prove that $R[M]$ satisfies the D-property if and only if $M$ is a valuation monoid and coprime nonzero elements of $R[M]$ remain coprime in $K[M]$. This includes the polynomial case and reduces over fields to the valuation property of $M$. We also prove that every TIDF D-domain is a PSP-domain and that every polynomial ring over it is TIDF. We exhibit a TIDF-domain that is not pre-Schreier although its polynomial ring satisfies the D-property. Finally, we construct a reduced antimatter D-monoid with a three-element subset having infinitely many maximal common divisors. This monoid cannot be realized as the reduced multiplicative monoid of a domain. Whether the IDF property ascends to polynomial rings over IDF D-domains remains open.
	\end{abstract}

	\maketitle

	\bigskip
	\section{Introduction}
	\label{sec:intro}
	
	Euclid's lemma connects divisibility of a product with coprimality of its factors. Its usual formulation for the integers extends to an integral domain $R$ as follows: if $a$, $b$, and $c$ in $R \setminus \{0\}$ satisfy $a \mid_R bc$ and every common divisor of $a$ and $b$ is a unit then $a \mid_R c$. We say that $R$ satisfies the \emph{D-property} when this implication holds. The same definition applies to commutative cancellative monoids, including the multiplicative monoids of nonzero elements of integral domains. A domain or monoid satisfying this property is called a \emph{D-domain} or a \emph{D-monoid}, respectively. The D-property does not require the existence of greatest common divisors or factorizations into atoms. Every atom of a D-monoid is prime, so every atomic D-monoid is a unique factorization monoid. Thus, the principal questions concerning the D-property arise in the absence of atomicity.
	\smallskip
	
	The D-property is equivalent to the two-generator form of the PSP property, also called the linear PSP property~\cite{HZ15}. A finite tuple of nonzero elements of a domain is \emph{primitive} if its entries have no common nonunit divisor, and it is \emph{superprimitive} if every element of the fraction field that sends all its entries into the domain already belongs to the domain. A \emph{PSP-domain} is a domain in which every primitive tuple is superprimitive. The D-property requires this condition only for pairs. Consequently, every maximal common divisor of a pair in a D-domain is a greatest common divisor, while the PSP property gives the corresponding conclusion for every finite nonempty tuple. In particular, a D-domain in which every pair has a maximal common divisor is a GCD-domain.
	\smallskip
	
	An integral domain is \emph{irreducible-divisor-finite}, or \emph{IDF}, if each nonzero element has only finitely many irreducible divisors up to associates. The IDF property does not in general ascend to polynomial rings. Eftekhari and Khorsandi~\cite{EK18} characterized this ascent by \emph{MCD-finiteness}: an IDF-domain $R$ has an IDF polynomial ring if and only if every finite nonempty tuple of nonzero elements of $R$ has finitely many maximal common divisors up to associates. Gonzalez and Panpaliya~\cite[Theorem~4.2]{GP25} showed that the IDF property can fail to ascend even among GL-domains. Whether it ascends to polynomial rings over IDF D-domains remains open. For further work on the IDF-ascent problem, see Malcolmson and Okoh~\cite{MO09}, Zafrullah~\cite{mZ17}, and Gotti and Zafrullah~\cite{GZ23}. Blitz, Gotti, Han, and Liang~\cite{BGHL26} studied the MCD and MCD-finite properties in monoid algebras, proving in particular that MCD-finiteness ascends to polynomial extensions.
	\smallskip
	
	We obtain a positive result under an additional factorization hypothesis. A domain is \emph{Furstenberg} if every nonzero nonunit has an irreducible divisor, and it is \emph{TIDF} if it is both Furstenberg and IDF. We prove that every TIDF D-domain is a PSP-domain and that every polynomial ring over such a domain is TIDF (Theorem~\ref{thm:tidf-d-is-psp}). The D-property makes every irreducible prime, allowing finite prime avoidance to produce a coprime pair in each nonzero finitely generated primitive ideal.
	\smallskip
	
	The examples in Section~\ref{sec:polynomial} illustrate the difference between ascent of the D-property and ascent of the IDF property. We construct antimatter PSP-domains $R$ for which $R[X]$ is IDF but does not satisfy the D-property (Example~\ref{ex:valuation-pullback} and Remark~\ref{ex:power-series-d-nonascent}). We also prove that the D-property of $R[X]$ forces $R$ to be integrally closed (Proposition~\ref{prop:polynomial-forces-normal}).
	\smallskip
	
	Related questions arise for factorization in monoid algebras. Coykendall and Gotti~\cite{CG19} constructed atomic monoids whose monoid algebras are nonatomic, and Gotti and Rabinovitz~\cite{GR25} obtained a rank-one example whose algebra over every integral domain is nonatomic. Gotti~\cite{fG22} studied factorization and closure properties of Puiseux algebras and constructed families of antimatter examples.
	\smallskip
	
	We characterize the Puiseux monoid algebras that satisfy the D-property. For a field $F$ and a Puiseux monoid $M$, we prove that $F[M]$ satisfies the D-property if and only if $M$ is a valuation monoid (Theorem~\ref{thm:puiseux-d-characterization}). More generally, for a domain $R$ with fraction field $K$ and a nontrivial Puiseux monoid $M$, the algebra $R[M]$ satisfies the D-property if and only if $M$ is a valuation monoid and every pair of coprime nonzero elements of $R[M]$ remains coprime in $K[M]$ (Theorem~\ref{thm:puiseux-coefficient-criterion}). The case $M = \nn_0$ gives the polynomial characterization in Corollary~\ref{cor:polynomial-d-characterization}.
	\smallskip
	
	If $R[M]$ satisfies the D-property for nontrivial $M$ then $R$ is a PSP-domain; if $R$ is also IDF then every polynomial ring over $R$ is IDF (Corollary~\ref{cor:d-algebra-coefficients}). We further show that if $R[X]$ satisfies the D-property then, for every Puiseux monoid $M$, the algebra $R[M]$ satisfies the D-property if and only if $M$ is a valuation monoid (Theorem~\ref{thm:d-polynomial-to-puiseux}). Example~\ref{ex:non-schreier-polynomial-ascent} gives a TIDF PSP-domain $R$ that is not pre-Schreier although $R[X]$ satisfies the D-property.
	\smallskip
	
	Finally, we show that the IDF and D-properties do not imply MCD-finiteness for monoids. In Section~\ref{sec:monoid}, we construct a reduced, antimatter (hence IDF) D-monoid in which some three-element subset has infinitely many MCDs (Theorem~\ref{thm:monoid}). This is the smallest possible cardinality, since a pair in a D-monoid has at most one MCD up to associates. The constructed monoid is not isomorphic to the reduced multiplicative monoid of an integral domain (Proposition~\ref{prop:nonrealizable}), so it does not settle the question for domains.

	\bigskip
	\section{Preliminaries}
	\label{sec:prelim}
	
	Throughout this paper, all monoids are commutative and cancellative, and all rings are commutative with identity. We use the term \emph{domain} for an integral domain. We write $\nn$ for the positive integers and set $\nn_0 \coloneqq \nn \cup \{0\}$. The symbols $\zz$, $\qq$, and $\rr$ denote the integers, rational numbers, and real numbers, respectively. Unless specified otherwise, $X$ denotes an indeterminate. For integers $a \le b$, we set $\ldb a, b \rdb \coloneqq [a, b] \cap \zz$. For $S \subseteq \rr$ and $r \in \rr$, we write $S_{\ge r} \coloneqq \{s \in S : s \ge r\}$.

	\medskip
	\subsection{Monoids}
	
	Let $M$ be a monoid, written multiplicatively with identity $1$. We denote its group of units by $M^\times$. For $a$ and $b$ in $M$, we say that $a$ \emph{divides} $b$, and write $a \mid_M b$, if $b = ac$ for some $c \in M$. The elements $a$ and $b$ are \emph{associates}, written $a \simeq_M b$, if $a = ub$ for some $u \in M^\times$. They are \emph{coprime}, written $a \perp_M b$, if every common divisor is a unit. We omit subscripts when the ambient monoid is clear. When $M$ is written additively, its identity is denoted by $0$, and $a \mid_M b$ means that $b = a + c$ for some $c \in M$.
	\smallskip
	
	An \emph{atom}, or \emph{irreducible element}, of $M$ is a nonunit $a \in M$ such that every equality $a = bc$ with $b$ and $c$ in $M$ forces $b$ or $c$ to be a unit. We denote the set of atoms of $M$ by $\mathcal{A}(M)$. The monoid $M$ is \emph{atomic} if every nonunit is a finite product of atoms. A nonunit $p \in M$ is \emph{prime} if $p \mid_M ab$ implies $p \mid_M a$ or $p \mid_M b$ for all $a$ and $b$ in $M$. An atomic monoid is a \emph{unique factorization monoid}, or \emph{UFM}, if every factorization into atoms is unique up to reordering and replacing atoms by associates. A monoid is \emph{IDF} if each element has finitely many irreducible divisors up to associates. It is \emph{antimatter} if it has no atoms, and it is \emph{Furstenberg} if every nonunit has an atom as a divisor. A monoid that is both IDF and Furstenberg is called a \emph{TIDF-monoid}.
	\smallskip
	
	A monoid is \emph{reduced} if its only unit is the identity. Its \emph{reduction} is the quotient by association. We write $\gp(M)$ for its \emph{Grothendieck group}, obtained by formally adjoining inverses to its elements, and identify $M$ with its canonical image in $\gp(M)$. A monoid is \emph{torsion-free} if $a^n = b^n$ with $a$ and $b$ in $M$ and $n \in \nn$ implies $a = b$. A torsion-free monoid has \emph{rank one} if its Grothendieck group is isomorphic to a nontrivial additive subgroup of $\qq$. The monoid $M$ is \emph{divisible} if, for every $a \in M$ and $n \in \nn$, there exists $b \in M$ such that $b^n = a$; in additive notation, this means $nb = a$.
	\smallskip
	
	A \emph{positive monoid} is an additive submonoid of $\rr_{\ge 0}$, and a \emph{Puiseux monoid} is an additive submonoid of $\qq_{\ge 0}$. A \emph{valuation monoid} is a monoid in which, for every pair of elements, one divides the other. A Puiseux monoid is a valuation monoid if and only if it is closed under nonnegative differences.
	\smallskip
	
	Let $n \in \nn$, and let $a_1, \ldots, a_n \in M$. The tuple $(a_1, \ldots, a_n)$ is \emph{primitive} if its entries have no common nonunit divisor. A common divisor $d \in M$ of the tuple is a \emph{maximal common divisor}, or \emph{MCD}, if $(a_1/d, \ldots, a_n/d)$ is primitive. Equivalently, every common divisor divisible by $d$ is associated to $d$. A \emph{greatest common divisor}, or \emph{GCD}, is a common divisor divisible by every common divisor. The monoid $M$ has the \emph{MCD property} if every finite nonempty tuple has an MCD. It is a \emph{GCD-monoid} if every finite nonempty tuple has a GCD. Equivalently, it suffices that every pair have a GCD (the \emph{$2$-GCD property}), since iterating pairwise GCDs gives a GCD of every finite nonempty tuple. The monoid $M$ is \emph{MCD-finite} if every finite nonempty tuple has finitely many MCDs up to associates. In a GCD-monoid every MCD of a tuple is associated to its GCD, so the monoid is MCD-finite. MCD-finiteness does not imply the MCD property: even in an MCD-finite domain, a pair can have no MCD~\cite[Example~3.3]{EK18}.
	\smallskip
	
	A finite nonempty tuple $(a_1, \ldots, a_n)$ in $M$ is \emph{superprimitive} if, for every $g \in \gp(M)$, the relations $ga_i \in M$ for all $i \in \ldb 1, n \rdb$ imply $g \in M$. A monoid is a \emph{PSP-monoid} if every primitive finite nonempty tuple is superprimitive. We denote the restriction of this condition to pairs by $\PSP_2$.

	\medskip
	\subsection{Integral Domains}
	
	Let $R$ be a domain. We denote its fraction field by $\Frac(R)$, its group of units by $R^\times$, and its multiplicative monoid of nonzero elements by $R^* \coloneqq R \setminus \{0\}$. The notions of divisibility, association, coprimality, atomicity, and unique factorization apply to $R$ through $R^*$. In particular, a \emph{unique factorization domain}, or \emph{UFD}, is a domain whose multiplicative monoid is a UFM.
	\smallskip
	
	A domain is IDF, antimatter, Furstenberg, or TIDF if its multiplicative monoid $R^*$ has the corresponding property. Every antimatter domain is IDF.
	\smallskip
	
	Primitivity, superprimitivity, and maximal and greatest common divisors of finite nonempty tuples of nonzero elements of $R$ are defined through $R^*$. Zero entries in a nonzero tuple may be discarded. A domain is an \emph{MCD-domain} if every such tuple has an MCD, and it is a \emph{$2$-MCD-domain} if every pair of nonzero elements has an MCD. The terms \emph{MCD-finite domain}, \emph{GCD-domain}, and \emph{PSP-domain} mean that $R^*$ has the corresponding property. The condition $\PSP_2$ for domains is also called the \emph{linear PSP property}~\cite{HZ15}. A \emph{B\'ezout domain} is a domain in which every finitely generated ideal is principal. A \emph{principal ideal domain}, or \emph{PID}, is a domain in which every ideal is principal. A \emph{valuation domain} is a domain whose multiplicative monoid of nonzero elements is a valuation monoid. The domain $R$ is \emph{integrally closed} if every element of $\Frac(R)$ that is a root of a monic polynomial in $R[X]$ belongs to $R$.
	\smallskip
	
	A nonzero polynomial in $R[X]$ is \emph{primitive} if its coefficients have no common nonunit divisor. The domain $R$ is a \emph{GL-domain} if the product of every two primitive polynomials in $R[X]$ is primitive~\cite{GP25}.
	\smallskip
	
	The ideal generated by a finite tuple is called primitive or superprimitive when the tuple has the corresponding property. A nonzero finitely generated ideal is primitive if and only if it is contained in no proper principal ideal. A nonzero $R$-submodule $I$ of $\Frac(R)$ is a \emph{fractional ideal} if $dI \subseteq R$ for some $d \in R^*$. For such an ideal, we write
	\[
	I^{-1} \coloneqq (R:I) = \{z \in \Frac(R) : zI \subseteq R\}.
	\]
	Thus, a nonzero finitely generated ideal $I$ of $R$ is superprimitive if and only if $I^{-1} = R$.
	\smallskip
	
	For a field $K$ and a totally ordered abelian group $\Gamma$ written additively, a \emph{valuation} is a group homomorphism $v \colon K^\times \to \Gamma$ satisfying $v(x + y) \ge \min\{v(x), v(y)\}$ for all $x$ and $y$ in $K^\times$ with $x + y \ne 0$. We set $v(0) = \infty$, and the image $v(K^\times)$ is the \emph{value group}.
	\smallskip
	
	The \emph{group of divisibility} of $R$ is $\Frac(R)^\times/R^\times$. For $x \in \Frac(R)^\times$, write $[x]$ for its class in this quotient. For $x$ and $y$ in $\Frac(R)^\times$, its order is defined by $[x] \le [y]$ if and only if $y/x \in R$. Its \emph{positive cone} is $R^*/R^\times$, the reduced multiplicative monoid of $R$.
	\smallskip
	
	We use the following criterion throughout the paper.
	\begin{theorem}[{Eftekhari--Khorsandi~\cite[Theorem~2.1]{EK18}}] \label{thm:ek}
		For an integral domain $R$, the following conditions are equivalent.
		\begin{enumerate}[(a)]
			\item The domain $R$ is IDF and MCD-finite.
			\smallskip
			
			\item The polynomial ring $R[X]$ is IDF.
			\smallskip
			
			\item Every polynomial ring over $R$ is IDF, including polynomial rings in infinitely many indeterminates.
		\end{enumerate}
	\end{theorem}
	
	Since every GCD-domain is MCD-finite, Theorem~\ref{thm:ek} implies that every polynomial ring over an IDF GCD-domain is IDF.
	\smallskip
	
	Let $R$ be a domain, and let $M$ be a torsion-free monoid written additively. The \emph{monoid algebra} $R[M]$ consists of formal sums $f = \sum_{m \in M} c_mX^m$, where $c_m \in R$ for every $m \in M$ and only finitely many coefficients are nonzero. Multiplication is determined by $X^mX^n = X^{m + n}$ for all $m$ and $n$ in $M$. Since $M$ is cancellative and torsion-free, the algebra $R[M]$ is a domain. The \emph{support} of $f$ is $\supp(f) \coloneqq \{m \in M : c_m \ne 0\}$. The \emph{augmentation map} is the homomorphism $\varepsilon \colon R[M] \rightarrow R$ defined as follows:
	\[
	\varepsilon\left( \sum_{m \in M}c_mX^m \right) = \sum_{m \in M} c_m.
	\]
	When $M$ is a positive monoid, we define the \emph{degree} of a nonzero element $f \in R[M]$ by $\deg(f) \coloneqq \max\supp(f)$. For nonzero elements $f$ and $g$ of $R[M]$, the identity $\deg(fg) = \deg(f) + \deg(g)$ implies that $R[M]^\times = R^\times$.

	\bigskip
	\section{The D-Property and Related Notions}
	\label{sec:d-related}
	
	This section relates the D-property to superprimitivity and maximal common divisors, establishes the PSP property for TIDF D-domains, and compares the D-property with the pre-Schreier property.
	
	\begin{definition}
		A monoid $M$ satisfies the \emph{D-property} if, for all $a$, $b$, and $c$ in $M$, the conditions $a \mid_M bc$ and $a \perp_M b$ imply $a \mid_M c$. Such a monoid is called a \emph{D-monoid}. A domain $R$ is a \emph{D-domain} if its multiplicative monoid $R^*$ satisfies the D-property.
	\end{definition}
	
	The D-property can be expressed as superprimitivity of coprime pairs.
	
	\begin{prop} \label{prop:d-psp2}
		A monoid satisfies the D-property if and only if it satisfies $\PSP_2$.
	\end{prop}
	
	\begin{proof}
		Let $M$ be a monoid, and assume that $M$ satisfies the D-property. Let $a$ and $b$ in $M$ be coprime, and let $g \in \gp(M)$ satisfy $ga \in M$ and $gb \in M$. The equality $b(ga) = a(gb)$ gives $a \mid_M b(ga)$, so the D-property implies $a \mid_M ga$. Cancellation yields $g \in M$.
		
		Conversely, suppose that every primitive pair is superprimitive. Take $a$, $b$, and $c$ in $M$ such that $a \mid_M bc$ and $a \perp_M b$. Set $g \coloneqq c/a \in \gp(M)$. Then $ga = c \in M$ and $gb = bc/a \in M$. Since the pair $(a, b)$ is superprimitive, it follows that $g \in M$ and hence $a \mid_M c$.
	\end{proof}

	\medskip
	\subsection{The PSP Property}
	
	The next two lemmas relate the PSP and D-properties to the MCD-finiteness condition in Theorem~\ref{thm:ek}.
	
	\begin{lemma} \label{lem:psp-unique}
		Let $M$ be a PSP-monoid. The following assertions hold.
		\begin{enumerate}[(1)]
			\item Every MCD of a finite nonempty tuple in $M$ is a GCD.
			\smallskip
			\item The monoid $M$ is MCD-finite.
		\end{enumerate}
	\end{lemma}
	
	\begin{proof}
		Let $M$ be a PSP-monoid, let $n \in \nn$, and let $s_1, \ldots, s_n \in M$. Suppose that $d \in M$ is an MCD of $(s_1, \ldots, s_n)$, and let $e \in M$ be a common divisor of this tuple. The tuple $(s_1/d, \ldots, s_n/d)$ is primitive. Moreover, multiplication by $d/e \in \gp(M)$ sends each entry into $M$, since $(d/e)(s_i/d) = s_i/e \in M$ for every $i \in \ldb 1, n \rdb$. The PSP property gives $d/e \in M$, so $e \mid_M d$. Hence $d$ is a GCD. Since GCDs are unique up to associates, each finite nonempty tuple has at most one MCD up to associates.
	\end{proof}
	
	\begin{lemma} \label{lem:pair-mcd}
		In a D-monoid, every MCD of a pair of elements is a GCD.
	\end{lemma}
	
	\begin{proof}
		Let $M$ be a D-monoid. For $a$ and $b$ in $M$, let $d \in M$ be an MCD of the pair $(a, b)$, and let $e \in M$ be a common divisor of this pair. Write $a = du = er$ and $b = dv = es$ with $u$, $v$, $r$, and $s$ in $M$. Then $u \perp_M v$, and $us = vr$. The D-property gives $u \mid_M r$. Writing $r = uh$ for some $h \in M$ and cancelling in $du = euh$, we obtain $d = eh$. Thus, $e \mid_M d$.
	\end{proof}
	
	For D-monoids, existence of pairwise MCDs therefore gives the GCD property.
	
	\begin{corollary} \label{cor:two-mcd-gcd}
		For a monoid $M$, the following conditions are equivalent.
		\begin{enumerate}[(a)]
			\item The monoid $M$ satisfies the D-property, and every pair of elements has an MCD.
			\smallskip
			
			\item The monoid $M$ is a GCD-monoid.
		\end{enumerate}
	\end{corollary}
	
	\begin{proof}
		(a) $\Rightarrow$ (b). In a D-monoid, every MCD of a pair is a GCD by Lemma~\ref{lem:pair-mcd}. Iterating pairwise GCDs gives the GCD property.
		\smallskip
		
		(b) $\Rightarrow$ (a). Suppose that $M$ is a GCD-monoid, and let $a$, $b$, and $c$ in $M$ satisfy $a \mid_M bc$ and $a \perp_M b$. The GCD cancellation law gives
		\[
		\gcd(ac, bc) \simeq_M c\gcd(a, b) \simeq_M c.
		\]
		Since $a$ divides both $ac$ and $bc$, it divides their GCD and hence divides $c$. Thus, $M$ satisfies the D-property. Every GCD is an MCD, so every pair has an MCD.
	\end{proof}
	
	We record the polynomial ascent theorem for GCD-domains (see Anderson~\cite{dA00}).
	
	\begin{theorem} \label{thm:gcd-polynomial-ascent}
		If $R$ is a GCD-domain then every polynomial ring over $R$ is a GCD-domain and hence satisfies the D-property.
	\end{theorem}
	
	\begin{proof}
		Iterating the one-variable theorem gives the finite-variable case. For infinitely many indeterminates, two polynomials involve only finitely many variables, and degrees in the remaining variables show that their divisors belong to the same finite-variable polynomial ring. The D-property now follows from Corollary~\ref{cor:two-mcd-gcd}.
	\end{proof}
	
	\begin{remark}
		Lemma~\ref{lem:psp-unique} guarantees that each finite nonempty tuple in a PSP-monoid has at most one MCD up to associates. Lemma~\ref{lem:pair-mcd} gives the same uniqueness conclusion for pairs in a D-monoid. Neither conclusion asserts the existence of an MCD; that is a separate hypothesis in Corollary~\ref{cor:two-mcd-gcd}.
	\end{remark}

	\medskip
	\subsection{TIDF and Atomic D-Domains}
	\label{sec:tidf}
	
	In a D-domain, the Furstenberg and IDF hypotheses allow us to find a coprime pair in every nonzero finitely generated primitive ideal. This reduction yields the PSP property.
	
	\begin{theorem} \label{thm:tidf-d-is-psp}
		Let $R$ be a TIDF-domain satisfying the D-property. The following assertions hold.
		\begin{enumerate}[(1)]
			\item The domain $R$ is a PSP-domain.
			\smallskip
			\item The domain $R$ is MCD-finite.
			\smallskip
			\item Every polynomial ring over $R$ is a TIDF-domain.
		\end{enumerate}
	\end{theorem}
	
	\begin{proof}
		Let $R$ be a TIDF-domain with the D-property. First, every irreducible of $R$ is prime. Indeed, let $p \in R$ be irreducible, and suppose that $p \mid_R bc$ for nonzero $b$ and $c$ in $R$. If $p \nmid_R b$ then $p \perp_R b$, because every nonunit divisor of $p$ is associate to $p$. The D-property gives $p \mid_R c$.
		
		Let $I$ be a nonzero finitely generated primitive ideal of $R$, and choose $a \in I \setminus \{0\}$. If $a$ is a unit then $I = R$, so $I$ is superprimitive. Suppose that $a$ is a nonunit. Since $R$ is TIDF, choose $r \in \nn$ and irreducibles $p_1, \ldots, p_r \in R$, one from each associate class of irreducible divisors of $a$. For every $j \in \ldb 1, r \rdb$, the ideal $p_jR$ is prime. Primitivity gives $I \nsubseteq p_jR$ for every $j \in \ldb 1, r \rdb$, so finite prime avoidance yields an element $b \in I$ such that
		\[
		b \notin \bigcup_{j = 1}^{r} p_jR.
		\]
		Note that $b \neq 0$ because $0 \in p_1R$. We claim that $a \perp_R b$. Otherwise, a common nonunit divisor $h \in R$ of $a$ and $b$ would have an irreducible divisor $p \in R$, because $R$ is Furstenberg. Then $p$ would divide both $a$ and $b$, and hence would be associate to some $p_j$ with $j \in \ldb 1, r \rdb$, contradicting the choice of $b$. This proves the claim. Proposition~\ref{prop:d-psp2} therefore gives
		\[
		R \subseteq I^{-1} \subseteq (a, b)^{-1} = R.
		\]
		Thus, $I$ is superprimitive, and $R$ is a PSP-domain.
		
		Lemma~\ref{lem:psp-unique} implies that $R$ is MCD-finite. The IDF property ascends to polynomial rings by Theorem~\ref{thm:ek}, and the Furstenberg property ascends by Gotti and Zafrullah~\cite[Proposition~4.7]{GZ23}. Hence every polynomial ring over $R$ is TIDF.
	\end{proof}
	
	The first step of the preceding proof does not require the IDF or Furstenberg hypotheses.
	
	\begin{corollary} \label{cor:atomic-d-ufd}
		Every atomic domain with the D-property is a UFD.
	\end{corollary}
	
	\begin{proof}
		Every irreducible is prime by the first paragraph of the proof of Theorem~\ref{thm:tidf-d-is-psp}. Atomicity now gives unique factorization.
	\end{proof}

	\medskip
	\subsection{The Pre-Schreier Property}
	
	Following Cohn~\cite{pC68}, an element $a$ of a monoid $M$ is called \emph{primal} if, whenever $a \mid_M bc$ for $b$ and $c$ in $M$, there exist $b'$ and $c'$ in $M$ such that
	\[
	a = b'c'
	\]
	with $b' \mid_M b$ and $c' \mid_M c$. A monoid is \emph{pre-Schreier} if every nonunit is primal. An integral domain is \emph{pre-Schreier} if its multiplicative monoid of nonzero elements is pre-Schreier. We use the terminology of Zafrullah~\cite{mZ87}; the underlying class was considered by Cohn~\cite{pC68}.
	
	Primality gives the divisibility implication required by the D-property.
	
	\begin{prop} \label{prop:pre-schreier-implies-d}
		Every pre-Schreier monoid satisfies the D-property.
	\end{prop}
	
	\begin{proof}
		Let $M$ be a pre-Schreier monoid, and take $a$, $b$, and $c$ in $M$ such that $a \mid_M bc$ and $a \perp_M b$. If $a$ is a unit then $a \mid_M c$. Suppose that $a$ is a nonunit. As $a$ is primal, there exist $b'$ and $c'$ in $M$ such that
		\[
		a = b'c'
		\]
		with $b' \mid_M b$ and $c' \mid_M c$. The element $b'$ divides both $a$ and $b$, so $b'$ is a unit. Hence $a \simeq_M c'$, and therefore $a \mid_M c$.
	\end{proof}
	
	The converse of Proposition~\ref{prop:pre-schreier-implies-d} fails even for Puiseux monoids.
	
	\begin{example} \label{ex:d-not-pre-schreier}
		Consider the Puiseux monoid
		\[
		M \coloneqq \zz[1/2]_{\ge 0} \cup \zz[1/2, 1/3]_{\ge 1}.
		\]
		This monoid satisfies the D-property but is not pre-Schreier.
	\end{example}
	
	\begin{proof}
		
		Call a rational number \emph{dyadic} if it belongs to $\zz[1/2]$. The union defining $M$ is closed under addition: a sum of dyadic elements is dyadic, and a sum involving an element of the second set belongs to that set.
		
		Every nonzero $q \in M$ is divisible in $M$ by $2^{-n}$ for all sufficiently large $n \in \nn$. For a positive dyadic $q$, choose $n$ such that $2^{-n} \le q$. If $q$ is not dyadic then $q > 1$, and we choose $n$ such that $2^{-n} \le q - 1$. In either case, $q - 2^{-n} \in M$. Thus, every finite nonempty tuple of nonzero elements of $M$ has a nonzero common divisor.
		
		It follows that every primitive finite nonempty tuple in $M$ contains the unit $0$ and is therefore superprimitive. Hence $M$ is a PSP-monoid and satisfies the D-property by Proposition~\ref{prop:d-psp2}.
		
		To prove that $M$ is not pre-Schreier, observe that
		\[
		1 \mid_M \frac43 + \frac43
		\]
		because $8/3 - 1 = 5/3 \in M$. If $1$ were primal then there would exist $u$ and $v$ in $M$ such that
		\[
		1 = u + v
		\]
		with $u \mid_M 4/3$ and $v \mid_M 4/3$. Since $0 \le u \le 1$ and $0 \le v \le 1$, both $u$ and $v$ are dyadic. If $u > 1/3$ then $4/3 - u$ is smaller than $1$ and is not dyadic, contradicting $4/3 - u \in M$. Moreover, $u \ne 1/3$ because $1/3 \notin M$. Therefore $u < 1/3$, and the same argument gives $v < 1/3$. This contradicts $u + v = 1$, proving that $1$ is not primal.
	\end{proof}
	
	An integral domain is a \emph{Schreier domain} if it is integrally closed and pre-Schreier. Cohn's polynomial ascent theorem provides a sufficient condition for the D-property that does not require atomicity~\cite[Theorem~2.7]{pC68}.
	
	\begin{prop} \label{prop:schreier-polynomial-ascent}
		The following assertions hold.
		\begin{enumerate}[(1)]
			\item Every $2$-MCD-domain satisfying the D-property is a Schreier domain.
			\smallskip
			
			\item Every polynomial ring over a Schreier domain is Schreier and hence satisfies the D-property.
		\end{enumerate}
	\end{prop}
	
	\begin{proof}
		(1) Let $R$ be a $2$-MCD-domain satisfying the D-property. Fix $p$, $a$, and $b$ in $R^*$ with $p \mid_R ab$, and choose an MCD $d \in R^*$ of $p$ and $a$. Write $p = dp'$ and $a = da'$ with $p'$ and $a'$ in $R^*$. Since $d$ is an MCD, we have $p' \perp_R a'$. Cancelling $d$ gives $p' \mid_R a'b$, and the D-property gives $p' \mid_R b$. Thus,
		\[
		p = dp'
		\]
		with $d \mid_R a$ and $p' \mid_R b$. Hence $p$ is primal, proving that $R$ is pre-Schreier. By Corollary~\ref{cor:two-mcd-gcd}, the domain $R$ is a GCD-domain. Every GCD-domain is integrally closed~\cite{dA00}, so $R$ is Schreier.
		\smallskip
		
		(2) Let $R$ be a Schreier domain. Cohn's theorem states that $R[X]$ is Schreier whenever $R$ is Schreier. Iteration proves the assertion for finitely many indeterminates. For an arbitrary set of indeterminates, each divisibility identity is contained in a polynomial ring in finitely many of them, where the pre-Schreier property gives the required decomposition. Likewise, an element of the fraction field and an equation witnessing its integrality involve only finitely many indeterminates, so integral closedness reduces to the finite-variable case. Thus, every polynomial ring over $R$ is Schreier. Its D-property follows from Proposition~\ref{prop:pre-schreier-implies-d}.
	\end{proof}

	\bigskip
	\section{Polynomial Extensions: Examples and Necessary Conditions}
	\label{sec:polynomial}
	
	We examine polynomial extensions through examples and a necessary condition before establishing the general characterization for Puiseux monoid algebras in Section~\ref{sec:monoid-algebras}.
	
	Proposition~\ref{prop:schreier-polynomial-ascent} yields polynomial ascent of the D-property for D-domains in which every pair of nonzero elements has an MCD. Without this existence hypothesis, ascent can fail even for PSP-domains. The following example also separates ascent of the D-property from ascent of the IDF property.
	
	\begin{example} \label{ex:valuation-pullback}
		There exists an antimatter PSP-domain $R$ such that $R[X]$ is IDF but does not satisfy the D-property.
	\end{example}
	
	\begin{proof}
		Let $K \coloneqq \qq(\sqrt2)$, fix compatible rational powers of an indeterminate $t$, and set
		\[
		V \coloneqq \bigcup_{n \in \nn} K[t^{1/n}]_{(t^{1/n})},
		\qquad \mathfrak m \coloneqq \{f \in V : v(f) > 0\},
		\qquad \text{and} \qquad R \coloneqq \qq + \mathfrak m,
		\]
		where $v \colon \Frac(V) \to \qq \cup \{\infty\}$ is the $t$-adic valuation normalized by $v(t) = 1$, with $v(0) = \infty$. The family defining $V$ is directed by divisibility of the indices. Every finite collection of elements of its fraction field lies in $K(t^{1/n})$ for some $n \in \nn$, where membership in $K[t^{1/n}]_{(t^{1/n})}$ is determined by nonnegative valuation. Consequently, $V$ is a valuation domain with value group $\qq$, maximal ideal $\mathfrak m$, and residue field $K$.
		
		An element of $R$ is a unit precisely when its residue is nonzero. Indeed, the inverse in $V$ of an element with nonzero rational residue again has rational residue and hence belongs to $R$. Every finite tuple of nonzero nonunits of $R$ has a common nonunit divisor $t^q$, where $q \in \qq_{> 0}$ is smaller than the valuation of each entry. Thus, every primitive tuple contains a unit, and therefore $R$ is a PSP-domain. In particular, $R$ satisfies the D-property.
		
		Each nonzero nonunit $f \in R$ factors as $t^q(f/t^q)$ with both factors nonunits for every $q \in \qq$ with $0 < q < v(f)$. Hence $R$ is antimatter and therefore IDF. By Lemma~\ref{lem:psp-unique}, it is also MCD-finite. Theorem~\ref{thm:ek} now implies that $R[X]$ is IDF.
		
		In $R[X]$, consider
		\[
		A \coloneqq X^2 - 2,
		\qquad B \coloneqq tX - \sqrt2t,
		\qquad \text{and} \qquad C \coloneqq tX + \sqrt2t.
		\]
		These polynomials satisfy $BC = t^2A$. A constant common divisor of $A$ and $B$ is a unit because $A$ is monic. A nonconstant common divisor must be linear, say $rX + s$ with $r \in R^*$ and $s \in R$. Since it divides the monic polynomial $A$, its leading coefficient $r$ is a unit of $R$. Comparing coefficients in $B = w(rX + s)$, where $w \in R$, yields $s/r = -\sqrt2$. This contradicts $\sqrt2 \notin R$. Thus, $A \perp_{R[X]} B$, whereas $A \nmid_{R[X]} C$ by degree. Therefore $R[X]$ does not satisfy the D-property.
		
		Although $A$ and $B$ are coprime in $R[X]$, they have the common nonconstant divisor $X - \sqrt{2}$ in $\Frac(R)[X]$, since $\sqrt{2} = (\sqrt{2}t)/t \in \Frac(R)$. This failure of coprimality to persist over the fraction field motivates the condition in Corollary~\ref{cor:polynomial-d-characterization}.
	\end{proof}
	
	The same construction has a formal power series realization.
	
	\begin{remark} \label{ex:power-series-d-nonascent}
		Let $K \coloneqq \qq(\sqrt{2})$, let $t$ be an indeterminate, and choose compatible roots $t^{1/n}$ for each $n \in \nn$. Set
		\[
		P \coloneqq \bigcup_{n \in \nn} K[\![t^{1/n}]\!]
		\qquad \text{and} \qquad
		R \coloneqq \{f \in P : f(0) \in \qq\}.
		\]
		The power series rings defining $P$ form a directed family, and the units of $R$ are precisely the series with nonzero constant term. Indeed, the inverse of such a series belongs to the same power series ring and again has rational constant term. Let $k \in \nn$, and let $f_1, \ldots, f_k$ be nonzero nonunits of $R$. Choose $N \in \nn$ such that every $f_i$ with $i \in \ldb 1, k \rdb$ belongs to $K[\![t^{1/N}]\!]$. Then $t^{1/(2N)}$ divides $f_i$ in $R$ for every $i \in \ldb 1, k \rdb$, and each quotient has zero constant term. Thus, every primitive tuple contains a unit, and each nonzero nonunit factors into two nonunits. Therefore $R$ is an antimatter PSP-domain, and $R[X]$ is IDF by Lemma~\ref{lem:psp-unique} and Theorem~\ref{thm:ek}. Both $t$ and $\sqrt{2}t$ belong to $R$, whereas $\sqrt{2} \notin R$, so the same polynomials and coefficient comparison as in Example~\ref{ex:valuation-pullback} show that $R[X]$ does not satisfy the D-property.
	\end{remark}
	
	The valuation domain used in the construction gives a contrasting example of polynomial ascent.
	
	\begin{remark} \label{rem:antimatter-valuation-polynomial-ascent}
		Atomicity is not needed for positive polynomial ascent. The valuation domain $V$ in Example~\ref{ex:valuation-pullback} is a GCD-domain, so every polynomial ring over $V$ satisfies the D-property by Theorem~\ref{thm:gcd-polynomial-ascent}. Nevertheless, $V$ is antimatter. Indeed, for a nonzero nonunit $f \in V$, choose $q \in \qq$ with $0 < q < v(f)$. The factorization $f = t^q(t^{-q}f)$ has two nonunit factors in $V$.
	\end{remark}
	
	The coefficient domains in Example~\ref{ex:valuation-pullback} and Remark~\ref{ex:power-series-d-nonascent} are not integrally closed. The next proposition shows that integral closedness is necessary for polynomial ascent of the D-property.
	
	\begin{prop} \label{prop:polynomial-forces-normal}
		Let $R$ be a domain. If $R[X]$ satisfies the D-property then $R$ is integrally closed.
	\end{prop}
	
	\begin{proof}
		Set $K \coloneqq \Frac(R)$. Suppose, by way of contradiction, that $\alpha = u/v \in K$ is integral over $R$ but does not belong to $R$, where $u$ and $v$ belong to $R^*$. Choose a monic polynomial $f \in R[X]$ of degree $n \in \nn$ such that $f(\alpha) = 0$, and set $L \coloneqq vX - u \in R[X]$. Division by $X - \alpha$ gives $f = (X - \alpha)q$, where $q \in R[\alpha][X]$. Synthetic division expresses each coefficient of $q$ as an $R$-linear combination of $1, \alpha, \ldots, \alpha^{n - 1}$. Hence $v^{n - 1}q \in R[X]$, and $v^n f = L(v^{n - 1}q)$, so $L \mid_{R[X]} v^n f$.
		
		Every constant common divisor of $L$ and $f$ is a unit because $f$ is monic. A nonconstant common divisor must be linear, say $g = aX + b \in R[X]$ with $a \in R^*$ and $b \in R$. Its leading coefficient $a$ is a unit because $g \mid_{R[X]} f$. Since $g$ and $L$ are associates over $K$, they have the same root, and therefore $\alpha = -b/a \in R$. This contradiction shows that $L \perp_{R[X]} f$. The D-property would then imply $L \mid_{R[X]} v^n$, contrary to the degree of $L$.
	\end{proof}
	
	The exact polynomial criterion is given in Corollary~\ref{cor:polynomial-d-characterization}, as a specialization of Theorem~\ref{thm:puiseux-coefficient-criterion}.

	\bigskip
	\section{Ascent to Monoid Algebras}
	\label{sec:monoid-algebras}
	
	We first consider atomic exponent monoids. The argument for Corollary~\ref{cor:atomic-d-ufd} applies to monoids as follows.
	
	\begin{prop} \label{prop:atomic-d-implies-uf}
		Let $M$ be an atomic monoid satisfying the D-property. Then $M$ is a unique factorization monoid.
	\end{prop}
	
	\begin{proof}
		Every atom of $M$ is prime. Indeed, if an atom $a \in M$ does not divide an element $b \in M$ then $a$ and $b$ have no common nonunit divisor. Whenever $a$ divides $bc$ for some $c \in M$, the D-property therefore forces $a$ to divide $c$.
		
		Since $M$ is atomic, it has unique factorization.
	\end{proof}
	
	The preceding proposition reduces the atomic case to polynomial algebras.
	\begin{prop} \label{prop:atomic-d-monoid-algebra}
		Let $M$ be a reduced torsion-free atomic monoid satisfying the D-property. The following assertions hold.
		\begin{enumerate}[(1)]
			\item For every GCD-domain $R$, the monoid algebra $R[M]$ is a GCD-domain.
			\smallskip
			\item For every field $F$, the monoid algebra $F[M]$ satisfies the D-property.
		\end{enumerate}
	\end{prop}
	
	\begin{proof}
		By Proposition~\ref{prop:atomic-d-implies-uf}, $M$ has unique factorization. As $M$ is also reduced, it is a free commutative monoid on its atoms; that is, $M \cong \bigoplus_{\alpha \in \mathcal{A}(M)} \nn_0$. Consequently,
		\[
		R[M] \cong R[X_\alpha : \alpha \in \mathcal{A}(M)],
		\]
		where the indeterminates are indexed by the atoms of $M$. Theorem~\ref{thm:gcd-polynomial-ascent} therefore gives both assertions.
	\end{proof}
	
	An \emph{affine monoid} is a finitely generated submonoid of $\zz^r$ for a nonnegative integer $r$. Every reduced affine monoid is atomic. Let $M$ be a reduced affine monoid satisfying the D-property. Then $M$ is isomorphic to $\nn_0^d$ for some $d \in \nn_0$, and, for every field $F$, its monoid algebra satisfies
	\[
	F[M] \cong F[X_1, \ldots, X_d].
	\]
	If $d \ge 2$ then distinct standard basis vectors of $\nn_0^d$ are incomparable under divisibility. Thus, the exponent monoids in these examples need not be valuation monoids.
	
	For Puiseux monoid algebras over fields, we obtain a characterization without an atomicity hypothesis. The next theorem shows that the D-property of the monoid algebra is equivalent to its B\'ezout property and to the valuation property of the exponent monoid.
	
	\begin{theorem} \label{thm:puiseux-d-characterization}
		Let $F$ be a field, and let $M$ be a Puiseux monoid. The following conditions are equivalent.
		
		\begin{enumerate}
			\item[(a)] The monoid $M$ is a valuation monoid.
			\smallskip
			
			\item[(b)] The monoid algebra $F[M]$ is a B\'ezout domain.
			\smallskip
			
			\item[(c)] The monoid algebra $F[M]$ satisfies the D-property.
		\end{enumerate}
	\end{theorem}
	
	Before proving Theorem~\ref{thm:puiseux-d-characterization}, we establish four lemmas. The first allows us to place finitely many elements of a valuation Puiseux monoid algebra in a common polynomial subring; we will also use it for arbitrary coefficient domains.
	
	\begin{lemma} \label{lem:valuation-finite-support}
		Let $M$ be a nontrivial valuation Puiseux monoid. For every finite subset $E \subseteq M$, there exists $d \in M \setminus \{0\}$ such that $E \subseteq d\nn_0$.
	\end{lemma}
	
	\begin{proof}
		If $E \subseteq \{0\}$ then every $d \in M \setminus \{0\}$ has the required property. Suppose that $E$ contains a positive element. Choose $N \in \nn$ such that $NE \subseteq \nn_0$, and let $g \in \nn$ be the greatest common divisor of the positive integers in $NE$. Set $d \coloneqq g/N$. Every element of $E$ is a nonnegative integral multiple of $d$. The subtraction form of the Euclidean algorithm computes $g$ from the positive elements of $NE$. Dividing each intermediate term by $N$ computes $d$ from the positive elements of $E$ using nonnegative differences. Since a valuation Puiseux monoid is closed under such differences, every step remains in $M$. Hence $d \in M \setminus \{0\}$, as required.
	\end{proof}
	
	\begin{lemma} \label{lem:binomial-divisor-support}
		Let $F$ be a field, let $n$ be a positive integer, and let $f \in F[X]$ be a nonconstant divisor of $X^n - 1$. Then $n$ belongs to the additive monoid generated by $\supp(f)$.
	\end{lemma}
	
	\begin{proof}
		Set $c \coloneqq f(0) \in F$. Since $f$ divides $X^n - 1$, its constant term $c$ is nonzero. Define a polynomial $q \in F[X]$ by
		\[
		q \coloneqq 1 - c^{-1}f.
		\]
		The polynomial $q$ has zero constant term, satisfies $f = c(1 - q)$, and has support contained in $\supp(f)$. Write $X^n - 1 = fh$ for a polynomial $h \in F[X]$, and set $g \coloneqq - ch \in F[X]$. Then the identity $(1 - q)g = 1 - X^n$ holds. Comparing constant terms gives $g(0) = 1$, while degrees give $\deg(g) = n - \deg(f) < n$. Expanding in the formal power series ring $F[\![X]\!]$ yields
		\[
		\sum_{j = 0}^{\infty}gX^{jn} = \frac{g}{1 - X^n} = \frac{1}{1 - q} = \sum_{j = 0}^{\infty}q^j.
		\]
		The coefficient of $X^n$ on the left is $1$, since only the summand with index $j = 1$ contributes to it. Hence the coefficient of $X^n$ on the right is nonzero. Since $q(0) = 0$, only finitely many powers of $q$ contribute to this coefficient. Therefore there is a positive integer $j$ such that $X^n$ occurs in $q^j$. The exponent $n$ is a sum of $j$ elements of $\supp(q)$, so it belongs to the monoid generated by $\supp(f)$.
	\end{proof}
	
	\begin{lemma} \label{lem:consecutive-integers}
		Let $M$ be a Puiseux monoid. If $a$ and $a + 1$ belong to $M$ for some nonnegative rational number $a$ then $M$ contains two consecutive positive integers.
	\end{lemma}
	
	\begin{proof}
		If $a = 0$ then $M$ contains both $1$ and $2$.
		
		Suppose now that $a > 0$, and write $a = u/v$ for positive integers $u$ and $v$. Closure under addition gives
		\[
		va = u \in M
		\qquad \text{and} \qquad
		(a + 1) + (v - 1)a = va + 1 = u + 1 \in M.
		\]
		Thus, the monoid $M$ contains the consecutive positive integers $u$ and $u + 1$.
	\end{proof}
	
	\begin{lemma} \label{lem:binomial-quotients}
		Let $F$ be a field, let $n$ be a positive integer, and let $M$ be a Puiseux monoid containing $n$ and $n + 1$. If two primes $p$ and $q$ satisfy
		\[
		p > q > n^2 - n - 1,
		\]
		then the four polynomials
		\[
		X^p - 1, \qquad X^q - 1, \qquad \frac{X^{pq} - 1}{X^p - 1}, \qquad \text{and} \qquad
		\frac{X^{pq} - 1}{X^q - 1}
		\]
		belong to $F[M]$.
	\end{lemma}
	
	\begin{proof}
		It suffices to show that the supports of all four polynomials are contained within $M$.
		
		Let $r$ be an integer satisfying $r \ge n(n - 1)$. By the division algorithm, we can write $r = kn + s$ for nonnegative integers $k$ and $s$, with $s < n$. Then the inequality $k \ge n - 1 \ge s$ gives $k - s \ge 0$, hence
		\[
		r = (k - s)n + s(n + 1) \in M.
		\]
		Thus, every integer greater than $n^2 - n - 1 = n(n - 1) - 1$ belongs to $M$.
		
		In particular, both $p$ and $q$ belong to $M$. Moreover,
		\[
		\frac{X^{pq} - 1}{X^p - 1} = 1 + X^p + X^{2p} + \cdots + X^{(q - 1)p}
		\]
		and
		\[
		\frac{X^{pq} - 1}{X^q - 1} = 1 + X^q + X^{2q} + \cdots + X^{(p - 1)q}.
		\]
		Each exponent in these sums is a nonnegative integral multiple of an element of $M$, so it belongs to~$M$. The supports of all four polynomials are therefore contained in $M$.
	\end{proof}
	
	We can now prove the field case of the characterization.
	
	\begin{proof}[Proof of Theorem~\ref{thm:puiseux-d-characterization}]
		(a) $\Rightarrow$ (b). Suppose that $M$ is a valuation monoid. We show that every finitely generated ideal of $F[M]$ is principal. The zero ideal is principal, so it suffices to consider an ideal generated by nonzero elements. Fix a positive integer $k$, choose nonzero generators $f_1, \ldots, f_k \in F[M]$, and define the finite set $E \subseteq M$ by
		\[
		E \coloneqq \bigcup_{i = 1}^k \supp(f_i).
		\]
		If $E = \{0\}$ then each generator is a nonzero constant, and the ideal is $F[M]$. Suppose that $E \neq \{0\}$. Then $M$ is nontrivial, so Lemma~\ref{lem:valuation-finite-support} gives $d \in M \setminus \{0\}$ such that $E \subseteq d\nn_0$.
		
		It follows that $f_1, \ldots, f_k$ belong to the polynomial ring $F[X^d]$, which is a PID. Let $g \in F[X^d]$ be a greatest common divisor of $f_1, \ldots, f_k$ in $F[X^d]$. There exist $u_1, \ldots, u_k \in F[X^d]$ such that
		\[
		g = u_1f_1 + \cdots + u_kf_k,
		\]
		and $g$ divides $f_i$ in $F[X^d]$ for every $i \in \ldb 1, k \rdb$. As $F[X^d] \subseteq F[M]$, we conclude that
		\[
		f_1F[M] + \cdots + f_kF[M] = gF[M]
		\]
		is a principal ideal. Hence $F[M]$ is a B\'ezout domain, and condition~(a) implies condition~(b).
		\smallskip
		
		(b) $\Rightarrow$ (c). Suppose that $F[M]$ is a B\'ezout domain. We show that it satisfies the D-property. Let $A$, $B$, and $C$ in $F[M]$ be nonzero elements satisfying $A \mid_{F[M]} BC$, and suppose that $A$ and $B$ are coprime. By condition~(b), the ideal generated by $A$ and $B$ is principal. A generator of this ideal divides both $A$ and $B$, so it is a unit. Thus, there exist elements $U$ and $V$ in $F[M]$ such that
		\[
		UA + VB = 1.
		\]
		Multiplying this identity by $C$ gives $UAC + VBC = C$. Both summands on the left are divisible by $A$, so $A \mid_{F[M]} C$. Thus, the monoid algebra $F[M]$ satisfies the D-property.
		\smallskip
		
		(c) $\Rightarrow$ (a). We argue by contraposition. Suppose that $M$ is not a valuation monoid. We show that $F[M]$ does not satisfy the D-property. Choose elements $u$ and $v$ in $M$ such that $u < v$ but $u \nmid_M v$, that is, $v - u \notin M$. Scaling all exponents by $1/(v - u)$ induces an isomorphism of monoid algebras, so we may assume that $v - u = 1$. Thus, $u$ and $u + 1$ belong to $M$, while $1$ does not.
		
		By Lemma~\ref{lem:consecutive-integers}, there exists $n \in \nn$ such that $n \in M$ and $n + 1 \in M$. Choose primes $p$ and $q$ satisfying
		\[
		p > q > n^2 - n - 1
		\]
		such that $p \cdot 1_F \ne 0$. By Lemma~\ref{lem:binomial-quotients}, the elements
		\[
		A \coloneqq X^p - 1, \qquad B \coloneqq X^q - 1, \qquad \text{and} \qquad C \coloneqq
		\frac{X^{pq} - 1}{X^q - 1}
		\]
		belong to $F[M]$. Moreover,
		\[
		BC = X^{pq} - 1 = A\left(\frac{X^{pq} - 1}{X^p - 1}\right),
		\]
		and so $A \mid_{F[M]} BC$.
		
		We claim that every common divisor of $A$ and $B$ in $F[M]$ is a unit. Suppose, by way of contradiction, that $D \in F[M]$ is a nonunit common divisor. Write $A = DU$ and $B = DV$ for elements $U$ and $V$ in $F[M]$. Choose a positive integer $m$ divisible by the denominators of all exponents in the supports of $D$, $U$, and $V$. Scaling exponents defines an isomorphism $\varphi_m \colon F[M] \rightarrow F[mM]$ as follows:
		\[
		\varphi_m\left(\sum_{e \in E} c_eX^e\right) \coloneqq \sum_{e \in E} c_eX^{me},
		\]
		where $E \subseteq M$ is the finite support of the element being mapped and $c_e \in F$ for every $e \in E$. The images of $D$, $U$, and $V$ belong to $F[X]$. Set $D' \coloneqq \varphi_m(D) \in F[X]$. Then the polynomial $D'$ divides both $X^{mp} - 1$ and $X^{mq} - 1$ in $F[X]$. Since~$D$ is a nonunit, it is nonconstant, so $D'$ is nonconstant as well. The Euclidean algorithm for binomials gives
		\[
		\gcd_{F[X]}\big(X^{mp} - 1, X^{mq} - 1\big) = X^{\gcd(mp, mq)} - 1 = X^m - 1,
		\]
		and so $D'$ divides $X^m - 1$ in $F[X]$. It follows from Lemma~\ref{lem:binomial-divisor-support} that
		\[
		m = e_1 + \cdots + e_\ell
		\]
		for a positive integer $\ell$ and elements $e_1, \ldots, e_\ell \in \supp(D')$. For each index $i \in \ldb 1, \ell \rdb$, write $e_i = md_i$ with $d_i \in \supp(D)$. Dividing the displayed equality by $m$ gives
		\[
		1 = d_1 + \cdots + d_\ell \in M,
		\]
		contradicting our choice of $M$. This proves the claim.
		
		If $F[M]$ satisfied the D-property then $A$ would divide $C$. There would therefore exist an element $H \in F[M]$ satisfying $AH = C$. Applying the augmentation homomorphism and using $\varepsilon(A) = 0$ would give $\varepsilon(C) = 0$. However,
		\[
		C = 1 + X^q + X^{2q} + \cdots + X^{(p - 1)q},
		\]
		so $\varepsilon(C) = p \cdot 1_F \neq 0$, a contradiction. Hence $F[M]$ does not satisfy the D-property, proving the implication from condition~(c) to condition~(a) by contraposition.
	\end{proof}

	\medskip
	\subsection{Arbitrary Coefficient Domains}
	\label{subsec:further-d-ascent}
	
	We now extend the preceding discussion to Puiseux monoid algebras over arbitrary integral domains. Throughout this subsection, $R$ is an integral domain with fraction field $K$, and $M$ denotes a Puiseux monoid. For a nonzero element $f \in K[M]$, let $A_f$ denote the fractional ideal of $R$ generated by the coefficients of $f$. A nonzero element $f \in R[M]$ is \emph{primitive} if its coefficients have no common nonunit divisor in $R$.
	
	We first establish the divisibility consequence of the PSP property needed below.
	
	\begin{lemma} \label{lem:psp-multiplier}
		Suppose that $R$ is a PSP-domain. The following assertions hold.
		\begin{enumerate}[(1)]
			\item If $f \in R[X]$ is primitive and $g \in K[X]$ satisfies $fg \in R[X]$ then $g \in R[X]$.
			\smallskip
			\item For every Puiseux monoid $M$, if $f \in R[M]$ is primitive and $g \in K[M]$ satisfies $fg \in R[M]$ then $g \in R[M]$.
		\end{enumerate}
	\end{lemma}
	
	\begin{proof}
		We first prove assertion~(1). Let $f \in R[X]$ be primitive, and let $g \in K[X]$ satisfy $fg \in R[X]$. If $g = 0$ then the conclusion holds, so suppose that $g \ne 0$. Let $A \coloneqq A_f$ be the ideal of $R$ generated by the coefficients of $f$. Since $R$ is a PSP-domain and $f$ is primitive, every $z \in K$ satisfying $zA \subseteq R$ belongs to $R$. We claim that, for each $j \in \nn$, every $z \in K$ satisfying $zA^j \subseteq R$ also belongs to $R$. The case $j = 1$ is the preceding observation. Suppose that the claim holds for $j$ and that $zA^{j + 1} \subseteq R$. For each $a \in A$, the inclusion $(za)A^j \subseteq R$ and the induction hypothesis imply that $za \in R$. Thus, $zA \subseteq R$, and hence $z \in R$. This proves the claim by induction.
		
		Set $r \coloneqq \deg(g) \in \nn_0$. The Dedekind--Mertens content formula (see Anderson~\cite{dA00}), applied after clearing the denominators of $g$, yields
		\[
		A^{r + 1}A_g = A^rA_{fg} \subseteq R.
		\]
		Each coefficient $z \in K$ of $g$ therefore satisfies $zA^{r + 1} \subseteq R$. By the claim, every such coefficient belongs to $R$, proving the polynomial assertion.
		
		For assertion~(2), let $f \in R[M]$ be primitive and let $g \in K[M]$ satisfy $fg \in R[M]$. The case $g = 0$ is immediate, so assume that $g \ne 0$. Choose a positive integer $d$ clearing the denominators of all exponents in the supports of $f$ and $g$. Multiplying these exponents by $d$ turns $f$ and $g$ into ordinary polynomials, preserves their coefficients, and respects multiplication. The first assertion therefore shows that every coefficient of $g$ lies in $R$. Since $g \in K[M]$, this yields $g \in R[M]$.
	\end{proof}
	
	The following observation extends the necessary condition in Theorem~\ref{thm:puiseux-d-characterization} to arbitrary coefficient domains.
	
	\begin{lemma} \label{lem:d-puiseux-necessary-valuation}
		Let $M$ be a Puiseux monoid. If $R[M]$ satisfies the D-property then $M$ is a valuation monoid.
	\end{lemma}
	
	\begin{proof}
		Suppose that $M$ is not a valuation monoid. As in the proof of Theorem~\ref{thm:puiseux-d-characterization}, we may scale $M$ so that $1 \notin M$ while $n$ and $n + 1$ belong to $M$ for some $n \in \nn$. Choose primes $p > q > n^2 - n - 1$ with $p \cdot 1_R \ne 0$, and set
		\[
		A \coloneqq X^p - 1, \qquad B \coloneqq X^q - 1, \qquad \text{and} \qquad
		C \coloneqq \sum_{i = 0}^{p - 1}X^{iq}.
		\]
		All three elements lie in $R[M]$, and
		\[
		BC = A\sum_{j = 0}^{q - 1}X^{jp}.
		\]
		The proof of Theorem~\ref{thm:puiseux-d-characterization}, applied over $K$, shows that $A$ and $B$ have no nonconstant common divisor in $K[M]$. Thus, every common divisor of $A$ and $B$ in $R[M]$ is a constant. Since the constant coefficient of $A$ is $-1$, such a divisor is a unit of $R$. Hence $A \perp_{R[M]} B$.
		
		On the other hand, $A$ does not divide $C$ in $R[M]$. Indeed, evaluation at $X = 1$ sends $A$ to $0$ and $C$ to $p \cdot 1_R \ne 0$. These three elements therefore witness the failure of the D-property in $R[M]$.
	\end{proof}
	
	We next isolate the additional condition needed for a characterization.
	
	\begin{definition} \label{def:FC-property}
		Let $M$ be a nontrivial Puiseux monoid, and let $R$ be an integral domain with fraction field $K$. The pair $(R, M)$ has the \emph{FC-property} if coprime nonzero elements of $R[M]$ remain coprime in $K[M]$.
	\end{definition}
	
	We first record a consequence of the FC-property for the coefficient domain.
	
	\begin{lemma} \label{lem:FC-implies-PSP}
		Let $M$ be a nontrivial Puiseux monoid. If $(R, M)$ has the FC-property then $R$ is a PSP-domain.
	\end{lemma}
	
	\begin{proof}
		Suppose, by way of contradiction, that $R$ is not a PSP-domain. Then there exist $n \in \nn_0$, a primitive tuple $(a_0, \ldots, a_n)$ of nonzero elements of $R$, and an element $z \in K \setminus R$ such that $za_i \in R$ for every $i \in \ldb 0, n \rdb$. Fix $m \in M \setminus \{0\}$, and define $A$, $f$, and $g$ in $K[M]$ by
		\[
		A \coloneqq \sum_{i = 0}^n a_iX^{2im}, \qquad
		f \coloneqq (X^m - z)A, \qquad \text{and} \qquad g \coloneqq a_0(X^m - z).
		\]
		Both $f$ and $g$ belong to $R[M]$ and have the nonconstant common divisor $X^m - z$ in $K[M]$. The FC-property therefore gives a common nonunit divisor $q \in R[M]$ of $f$ and $g$. Define $b \in K[M]$ by
		\[
		b \coloneqq \frac{X^m - z}{q} = \frac{1}{a_0}\frac{g}{q} \in K[M].
		\]
		Then $Ab = f/q \in R[M]$.
		
		Let $c \in R^*$ be the leading coefficient of $q$. The leading coefficient of $b$ is $1/c$, and $\deg(b) \le m$. For $i \in \ldb 0, n \rdb$, the supports of the summands $a_iX^{2im}b$ are contained in the pairwise disjoint intervals $[2im, (2i + 1)m]$. The coefficient of $X^{2im + \deg(b)}$ in $Ab$ is therefore $a_i/c$ for every $i \in \ldb 0, n \rdb$. Since $Ab \in R[M]$, each of these coefficients belongs to $R$. As the $a_i$ have no common nonunit divisor, $c$ is a unit of $R$. Replacing $q$ by $c^{-1}q$ and $b$ by $cb$, we may assume that $q$ is monic. As a nonunit, $q$ is then nonconstant.
		
		Clearing the denominators of the exponents in the identity $bq = X^m - z$ gives
		\[
		B(X)Q(X) = X^N - z
		\]
		for some positive integer $N$, some polynomial $B \in K[X]$, and a monic nonconstant polynomial $Q \in R[X]$. The remainder on dividing $X^N$ by $Q$ in $K[X]$ is therefore $z$. On the other hand, since $Q$ is monic, division in $R[X]$ gives a remainder in $R[X]$. Uniqueness of the remainder yields $z \in R$, a contradiction. Hence $R$ is a PSP-domain.
	\end{proof}
	
	We can now establish the Puiseux ascent criterion in terms of coprimality.
	
	\begin{theorem} \label{thm:puiseux-coefficient-criterion}
		Let $R$ be an integral domain with fraction field $K$, and let $M$ be a nontrivial Puiseux monoid. The following conditions are equivalent.
		\begin{enumerate}[(a)]
			\item The monoid algebra $R[M]$ satisfies the D-property.
			\smallskip
			
			\item The monoid $M$ is a valuation monoid, and the pair $(R, M)$ has the FC-property.
		\end{enumerate}
	\end{theorem}
	
	\begin{proof}
		(a) $\Rightarrow$ (b). By Lemma~\ref{lem:d-puiseux-necessary-valuation}, $M$ is a valuation monoid. To prove the FC-property, let $f$ and $g$ in $R[M]$ be coprime nonzero elements, and let $q \in K[M]$ be a common divisor. Write $f = qu$ and $g = qv$ with $u$ and $v$ in $K[M]$. Choose $\alpha \in R^*$ such that $\alpha u \in R[M]$ and $\alpha v \in R[M]$. Then
		\[
		f(\alpha v) = g(\alpha u),
		\]
		so the D-property gives $f \mid_{R[M]} \alpha u$. Hence $\alpha/q \in R[M]$, and therefore $1/q \in K[M]$. Thus, $q$ is a unit of $K[M]$, proving that $f$ and $g$ remain coprime. Hence $(R, M)$ has the FC-property.
		\smallskip
		
		(b) $\Rightarrow$ (a). Take nonzero $a$, $b$, and $c$ in $R[M]$ such that $a \mid_{R[M]} bc$ and $a \perp_{R[M]} b$. Choose $t \in M$ with $t > \deg(a)$, and define $F \in R[M]$ by
		\[
		F \coloneqq a + X^t b.
		\]
		Such a $t$ exists because $M$ is nontrivial. The supports of $a$ and $X^t b$ are disjoint. Hence the nonzero coefficients of $F$ are precisely the nonzero coefficients of $a$ and $b$, with their multiplicities. A common nonunit divisor of these coefficients would divide both $a$ and $b$ in $R[M]$. Therefore $F$ is primitive.
		
		Since $(R, M)$ has the FC-property, $a \perp_{K[M]} b$. The implication from condition~(a) to condition~(c) in Theorem~\ref{thm:puiseux-d-characterization} shows that $K[M]$ satisfies the D-property. Hence $a \mid_{K[M]} c$, and we may write $c = aq$ with $q \in K[M]$. Also write $bc = ad$ with $d \in R[M]$. Cancellation gives $bq = d$, so
		\[
		Fq = (a + X^t b)q = c + X^t d \in R[M].
		\]
		By Lemma~\ref{lem:FC-implies-PSP}, $R$ is a PSP-domain. Since $F$ is primitive, Lemma~\ref{lem:psp-multiplier} yields $q \in R[M]$. Consequently, $a \mid_{R[M]} c$, as required.
	\end{proof}
	
	Taking $M = \nn_0$ gives the polynomial characterization anticipated in Section~\ref{sec:polynomial}.
	
	\begin{corollary} \label{cor:polynomial-d-characterization}
		Let $R$ be an integral domain with fraction field $K$. The following conditions are equivalent.
		\begin{enumerate}[(a)]
			\item The polynomial ring $R[X]$ satisfies the D-property.
			\smallskip
			
			\item Every pair of coprime nonzero polynomials in $R[X]$ remains coprime in $K[X]$.
		\end{enumerate}
	\end{corollary}
	
	\begin{proof}
		The monoid $\nn_0$ is a nontrivial valuation Puiseux monoid, and $R[\nn_0] = R[X]$. Therefore the equivalence follows from Theorem~\ref{thm:puiseux-coefficient-criterion}.
	\end{proof}
	
	The characterization also connects the D-property of a monoid algebra with IDF ascent from its coefficient domain.
	
	\begin{corollary} \label{cor:d-algebra-coefficients}
		Let $M$ be a nontrivial Puiseux monoid, and suppose that $R[M]$ satisfies the D-property. The following assertions hold.
		\begin{enumerate}[(1)]
			\item The domain $R$ is a PSP-domain.
			\smallskip
			\item The domain $R$ is MCD-finite.
			\smallskip
			\item If $R$ is IDF then every polynomial ring over $R$ is IDF.
		\end{enumerate}
	\end{corollary}
	
	\begin{proof}
		By Theorem~\ref{thm:puiseux-coefficient-criterion}, the pair $(R, M)$ has the FC-property. Lemma~\ref{lem:FC-implies-PSP} implies that $R$ is a PSP-domain, and Lemma~\ref{lem:psp-unique} gives MCD-finiteness. The final assertion follows from Theorem~\ref{thm:ek}.
	\end{proof}
	
	In particular, if $R$ is an IDF-domain for which $R[X]$ is not IDF then $R[X]$ does not satisfy the D-property.
	
	For the trivial Puiseux monoid $M = \{0\}$, we have $R[M] = R$, so the D-property of the monoid algebra is simply the D-property of $R$.
	
	For valuation Puiseux monoids, we obtain the following equivalent formulations of the FC-property.
	
	\begin{prop} \label{prop:FC-equivalence}
		Let $M$ be a nontrivial valuation Puiseux monoid. The following conditions are equivalent.
		\begin{enumerate}[(a)]
			\item The pair $(R, M)$ has the FC-property.
			\smallskip
			
			\item The domain $R$ is a PSP-domain, and whenever $f \in R[M]$ is primitive and a nonconstant $q \in K[M]$ divides $f$ in $K[M]$, there exists a primitive nonconstant $g \in R[M]$ that divides $q$ in $K[M]$.
			\smallskip
			
			\item The domain $R$ is a PSP-domain, and for every primitive $f \in R[X]$, every nonconstant $q \in K[X]$ dividing $f$ in $K[X]$, and every $m \in M \setminus \{0\}$, there exist $n \in \nn$ and a primitive nonconstant $h \in R[X]$ such that
			\[
			h(X) \mid_{K[X]} q(X^n)
			\]
			with $m/n \in M$.
		\end{enumerate}
	\end{prop}
	
	\begin{proof}
		(a) $\Rightarrow$ (b). By Lemma~\ref{lem:FC-implies-PSP}, $R$ is a PSP-domain. Let $f \in R[M]$ be primitive and let $q \in K[M]$ be a nonconstant divisor of $f$. Choose $\alpha \in R^*$ such that $\alpha q \in R[M]$. The nonconstant element $q$ divides both $f$ and $\alpha q$ in $K[M]$. By condition~(a), these two elements have a common nonunit divisor $g \in R[M]$. Since $g$ divides the primitive element $f$ in $R[M]$, it is primitive; as a nonunit, it is also nonconstant. Finally, $g \mid_{R[M]} \alpha q$ implies $g \mid_{K[M]} q$.
		\smallskip
		
		(b) $\Rightarrow$ (a). Let $a$ and $b$ in $R[M]$ be coprime nonzero elements. Choose $t \in M$ with $t > \deg(a)$. The supports of $a$ and $X^t b$ are disjoint, so the nonzero coefficients of $F \coloneqq a + X^t b$ are precisely the nonzero coefficients of $a$ and $b$, with their multiplicities. A common nonunit divisor of these coefficients would divide both $a$ and $b$. Therefore $F$ is primitive. If $a$ and $b$ had a nonconstant common divisor $q \in K[M]$ then $q$ would divide $F$. By condition~(b), there would be a primitive nonconstant $g \in R[M]$ dividing $q$ in $K[M]$. Since $R$ is a PSP-domain, Lemma~\ref{lem:psp-multiplier} would give $g \mid_{R[M]} a$ and $g \mid_{R[M]} b$, a contradiction. Since $K[M]^\times = K^\times$, this proves that $a$ and $b$ remain coprime in $K[M]$. Thus, $(R, M)$ has the FC-property.
		\smallskip
		
		(b) $\Rightarrow$ (c). Let $f \in R[X]$ be primitive, let $q \in K[X]$ be a nonconstant divisor of $f$, and let $m \in M \setminus \{0\}$. Since $m > 0$, the substitution $X \mapsto X^m$ preserves coefficients and gives
		\[
		q(X^m) \mid_{K[M]} f(X^m).
		\]
		By condition~(b), there is a primitive nonconstant $g \in R[M]$ such that $q(X^m) = gv$ for some $v \in K[M]$. By Lemma~\ref{lem:valuation-finite-support}, choose $d \in M \setminus \{0\}$ such that
		\[
		\{m\} \cup \supp(g) \cup \supp(v) \subseteq d \nn_0.
		\]
		Write $m = nd$ for some $n \in \nn$, and choose polynomials $h \in R[X]$ and $w \in K[X]$ such that $g = h(X^d)$ and $v = w(X^d)$. The substitution $X \mapsto X^d$ is injective, so
		\[
		q(X^n) = h(X)w(X).
		\]
		Moreover, $h$ is primitive and nonconstant, and $m/n = d \in M$. Together with the PSP hypothesis in (b), this proves (c).
		\smallskip
		
		(c) $\Rightarrow$ (b). The domain $R$ is a PSP-domain by (c). Let $f \in R[M]$ be primitive and suppose that $f = qv$ for a nonconstant $q \in K[M]$ and some $v \in K[M]$. By Lemma~\ref{lem:valuation-finite-support}, choose $m \in M \setminus \{0\}$ such that the supports of $f$, $q$, and $v$ are contained in $m\nn_0$. Choose a polynomial $F \in R[X]$ and polynomials $Q$ and $V$ in $K[X]$ such that
		\[
		f = F(X^m), \qquad q = Q(X^m), \qquad \text{and} \qquad v = V(X^m).
		\]
		Then $F$ is primitive, $Q$ is nonconstant, and $F = QV$. By condition~(c), there are $n \in \nn$ and a primitive nonconstant $h \in R[X]$ such that $m/n \in M$ and $h(X)$ divides $Q(X^n)$ in $K[X]$. Substituting $X^{m/n}$ for $X$ shows that $h(X^{m/n})$ is a primitive nonconstant element of $R[M]$ dividing $q$ in $K[M]$.
	\end{proof}
	
	When $R[X]$ satisfies the D-property, the valuation property of $M$ alone determines whether $R[M]$ satisfies the D-property.
	
	\begin{theorem} \label{thm:d-polynomial-to-puiseux}
		Suppose that $R[X]$ satisfies the D-property. For every Puiseux monoid $M$, the monoid algebra $R[M]$ satisfies the D-property if and only if $M$ is a valuation monoid.
	\end{theorem}
	
	\begin{proof}
		Necessity follows from Lemma~\ref{lem:d-puiseux-necessary-valuation}. Conversely, suppose that $M$ is a valuation monoid. If $M = \{0\}$ then $R[M] = R$, and the D-property descends from $R[X]$ to $R$ because every divisor of a nonzero constant polynomial is constant. Suppose that $M$ is nontrivial.
		
		Let $a$, $b$, $c$, and $d$ in $R[M] \setminus \{0\}$ satisfy $bc = ad$ and $a \perp_{R[M]} b$. Applying Lemma~\ref{lem:valuation-finite-support} to the union of their supports gives $q \in M \setminus \{0\}$ such that $a$, $b$, $c$, and $d$ belong to $R[X^q]$. Both $R[X^q]$ and $R[M]$ have unit group $R^\times$, so $a$ and $b$ are coprime in $R[X^q]$. Since $R[X^q]$ is isomorphic to $R[X]$, it satisfies the D-property. The identity $bc = ad$ therefore gives $a \mid_{R[X^q]} c$, and hence $a \mid_{R[M]} c$.
	\end{proof}
	
	Combining this theorem with Corollary~\ref{cor:polynomial-d-characterization} and Theorem~\ref{thm:puiseux-coefficient-criterion} shows that the FC-property passes from $(R, \nn_0)$ to $(R, M)$ for every nontrivial valuation Puiseux monoid $M$.
	
	The preceding theorem applies to coefficient domains outside the pre-Schreier class. The following example uses a pullback construction of the type considered by Anderson and Zafrullah~\cite[Example~2.11]{AZ07}.
	
	\begin{example} \label{ex:non-schreier-polynomial-ascent}
		Let $U$ and $T$ be indeterminates over $\qq$, and set
		\[
		A \coloneqq \zz_{(2)}, \qquad L \coloneqq \qq(U), \qquad V \coloneqq L[\![T]\!], \qquad \text{and} \qquad R \coloneqq A + TV.
		\]
		Here $A$ is the localization of $\zz$ at the prime ideal $(2)$. The domain $R$ is a TIDF PSP-domain that is not pre-Schreier, while $R[X]$ satisfies the D-property.
	\end{example}
	
	\begin{proof}
		The fraction field of $R$ is $K \coloneqq L((T))$. Indeed, the inclusion $R \subseteq V$ gives $\Frac(R) \subseteq \Frac(V)$, while $v = (Tv)/T \in \Frac(R)$ for every $v \in V$ gives the reverse inclusion.
		
		An element of $R$ is a unit precisely when its constant coefficient is a unit of $A$. Every nonunit belongs to $2R$, so every primitive finite tuple contains a unit. Thus, $R$ is a PSP-domain. Every nonzero element with nonzero constant coefficient is $2^n$ times a unit for some $n \in \nn_0$. Every nonzero element with zero constant coefficient remains a nonunit after division by $2$. Consequently, $2$ is the only atom up to associates, and it divides every nonzero nonunit. Hence $R$ is TIDF.
		
		To see that $R$ is not pre-Schreier, observe that $T \mid_R (UT)(T/U)$. Suppose that $T = rs$ for $r$ and $s$ in $R$ with $r \mid_R UT$ and $s \mid_R T/U$. The $T$-adic orders of $r$ and $s$ are nonnegative integers whose sum is $1$. If $r$ has order $0$ then $s \mid_R T/U$ gives $r/U \in R$, but its constant coefficient $r(0)/U$ does not belong to $A$. If $s$ has order $0$ then $r \mid_R UT$ gives $Us \in R$, but its constant coefficient $Us(0)$ does not belong to $A$. Both cases are impossible, so $T$ is not primal.
		
		We verify condition~(b) of Proposition~\ref{prop:FC-equivalence} for $M = \nn_0$. Let $f \in R[X]$ be primitive, and let $q \in K[X]$ be a nonconstant divisor of $f$. For a nonzero polynomial over $K$, the minimum $T$-adic order of its coefficients is additive under multiplication. To see this, rescale each factor by a power of $T$ so that its minimum coefficient order is zero. The reductions modulo $T$ are then nonzero polynomials over the field $L$, so their product is nonzero. This minimum is $0$ for $f$, because $f$ has a unit coefficient. Thus, multiplying $q$ by a suitable power of $T$ gives a polynomial $q_1 \in V[X]$ such that $h_1 \coloneqq f/q_1 \in K[X]$ also belongs to $V[X]$, and both have nonzero reductions modulo $T$. Writing a bar for this reduction, we obtain
		\[
		\overline{q_1}\, \overline{h_1} = \overline{f} \in A[X] \setminus \{0\}.
		\]
		Irreducible polynomials over $\qq$ remain irreducible over the purely transcendental extension $L = \qq(U)$. Unique factorization in $L[X]$ therefore implies that $\overline{q_1}$ is an $L^\times$-multiple of a polynomial in $\qq[X]$. Rescaling that polynomial by an integer power of $2$ gives a polynomial in $A[X]$ with a unit coefficient. Consequently, there exists $\lambda \in L^\times$ such that $g \coloneqq \lambda q_1 \in V[X]$ has reduction in $A[X]$ with a unit coefficient. Then $g \in R[X]$ is primitive and nonconstant, and it is associate to $q$ in $K[X]$. This proves condition~(b), so $(R, \nn_0)$ has the FC-property. Corollary~\ref{cor:polynomial-d-characterization} now shows that $R[X]$ satisfies the D-property.
		
		By Theorem~\ref{thm:d-polynomial-to-puiseux}, the monoid algebra $R[M]$ satisfies the D-property for every valuation Puiseux monoid $M$. In particular, this holds for $M = \zz[1/2]_{\ge 0}$.
	\end{proof}

	\bigskip
	\section{An Antimatter D-Monoid with Infinitely Many MCDs}
	\label{sec:monoid}
	
	By Lemma~\ref{lem:pair-mcd}, every pair in a D-monoid has at most one MCD up to associates. We show that a three-element subset can have infinitely many MCDs, even in a reduced antimatter monoid. Thus, three is the smallest possible cardinality for such a subset, and the IDF and D-properties do not imply MCD-finiteness for monoids. We also prove that the constructed monoid cannot be realized as the reduced multiplicative monoid of an integral domain.
	
	Label the coordinates of $\qq^7$ by $A$, $B$, $C$, $Z$, $U$, $V$, and $W$, in that order. Call $A$, $B$, and $C$ the \emph{base coordinates} and $Z$, $U$, $V$, and $W$ the \emph{auxiliary coordinates}. Let $\mathcal{F}$ be the family of nonempty unions of the six sets
	\[
	\{A\}, \quad \{B\}, \quad \{C\}, \quad
	\{A, B, Z, U\}, \quad \{A, C, Z, V\}, \quad \{B, C, Z, W\},
	\]
	and set
	\begin{equation} \label{eq:monoid}
		M \coloneqq \{0\} \cup \{v \in \qq_{\ge 0}^7 : \supp(v) \in \mathcal{F}\}.
	\end{equation}
	Here $\supp(v)$ denotes the set of coordinates at which $v$ is nonzero. More explicitly, $\mathcal{F}$ consists of the seven nonempty subsets of $\{A, B, C\}$ and the ten sets
	\[
	\begin{gathered}
		ABZU, \ ABCZU, \ ACZV, \ ABCZV, \ BCZW, \ ABCZW, \\
		ABCZUV, \ ABCZUW, \ ABCZVW, \ ABCZUVW,
	\end{gathered}
	\]
	where each string denotes the set of its listed coordinates. For a nonempty subset $J \subseteq \{U, V, W\}$, let $\kappa(J)$ be the union of the sets specified by
	\[
	\kappa(\{U\}) = \{A, B\}, \qquad
	\kappa(\{V\}) = \{A, C\}, \qquad \text{and} \qquad
	\kappa(\{W\}) = \{B, C\}.
	\]
	Then a subset $S \subseteq \{A, B, C, Z, U, V, W\}$ belongs to $\mathcal{F}$ if and only if either $S$ is a nonempty subset of $\{A, B, C\}$, or $Z \in S$, the set $J \coloneqq S \cap \{U, V, W\}$ is nonempty, and $\kappa(J) \subseteq S$.
	
	The base-coordinate sets $\{A, B\}$, $\{A, C\}$, and $\{B, C\}$ intersect pairwise but have empty total intersection. The auxiliary coordinates $U$, $V$, and $W$ record these pairs and control coprimality, while the common coordinate $Z$ provides the parameter for the infinite family of MCDs.
	
	\begin{theorem} \label{thm:monoid}
		The monoid $M$ in~\eqref{eq:monoid} has the following properties.
		\begin{enumerate}[(1)]
			\item The monoid $M$ is reduced, cancellative, and divisible, and it satisfies the D-property.
			\smallskip
			\item The monoid $M$ is antimatter and IDF.
			\smallskip
			\item Some three-element subset of $M$ has infinitely many MCDs.
			\smallskip
			\item The Grothendieck group of $M$ is $\qq^7$.
		\end{enumerate}
	\end{theorem}
	
	\begin{proof}
		Since $\mathcal{F}$ is closed under unions, the set $M$ is closed under addition. It is reduced and cancellative because it is a submonoid of $\qq_{\ge 0}^7$. Positive rational rescaling preserves support, so $M$ is divisible. In particular, every nonzero element $a \in M$ satisfies $a = a/2 + a/2$, and hence $M$ is antimatter and therefore IDF. Every vector in $\qq^7$ is the difference of two vectors with all coordinates positive. Thus, $\gp(M) = \qq^7$.
		
		To prove the D-property, observe that two nonzero elements $a$ and $b$ in $M$ have a nonzero common divisor if and only if their base supports intersect. Indeed, if both vectors are positive at a base coordinate then a sufficiently small positive rational multiple of its coordinate vector is a common divisor: subtracting it leaves both supports unchanged. Conversely, if $d \in M$ is a nonzero common divisor of $a$ and $b$ then $a - d$ and $b - d$ are coordinatewise nonnegative, so $a$ and $b$ are positive wherever $d$ is; since every nonzero element of $M$ is positive at some base coordinate, so is $d$. The same argument shows that a finite nonempty tuple of nonzero elements of $M$ has a nonzero common divisor if and only if there is a base coordinate at which every entry is positive.
		
		By Proposition~\ref{prop:d-psp2}, it suffices to prove that $q \in M$ whenever $a$ and $b$ are coprime elements of $M$ and $q \in \gp(M)$ satisfies $q + a \in M$ and $q + b \in M$. Fix such $a$, $b$, and $q$. If $a = 0$ or $b = 0$ then $q \in M$ directly. Suppose that $a$ and $b$ are nonzero. Their base supports are disjoint. An allowed support containing an auxiliary coordinate uses at least two base coordinates, so at most one of $a$ and $b$ has an auxiliary coordinate. Their full supports are therefore disjoint. Since the criterion is symmetric in $a$ and $b$, we may assume that $b$ is supported on the base coordinates.
		
		The disjointness of the supports of $a$ and $b$, together with the nonnegativity of $q + a$ and $q + b$, implies that $q$ is coordinatewise nonnegative. For a vector $v \in \qq_{\ge 0}^7$, let $P_v$ be its base support and let $J_v \coloneqq \supp(v) \cap \{U, V, W\}$. Since $b$ is supported on the base coordinates, $J_{q + b} = J_q$. If $J_q$ is empty then $q + b \in M$ forces the $Z$-coordinate of $q$ to be zero by the membership criterion for $\mathcal{F}$, and hence $q \in M$. If $J_q$ is nonempty then $q + b \in M$ forces the $Z$-coordinate of $q$ to be positive and gives
		\[
		\kappa(J_q) \subseteq P_q \cup P_b.
		\]
		Similarly, the condition $q + a \in M$ gives
		\[
		\kappa(J_q) \subseteq \kappa(J_q \cup J_a) \subseteq P_q \cup P_a.
		\]
		Since $P_a \cap P_b$ is empty, these inclusions imply $\kappa(J_q) \subseteq P_q$. Together with the positive $Z$-coordinate, the membership criterion for $\mathcal{F}$ yields $\supp(q) \in \mathcal{F}$, that is, $q \in M$. Thus, $M$ has the D-property.
		
		Finally, let $u_1$, $u_2$, and $u_3$ in $M$ be the indicator vectors of $ABZU$, $ACZV$, and $BCZW$, respectively, and put $a_i \coloneqq u_1 + u_i \in M$ for $i \in \ldb 1, 3 \rdb$. The supports of $a_1$, $a_2$, and $a_3$ are $ABZU$, $ABCZUV$, and $ABCZUW$, respectively, so these elements are distinct. For each rational number $t \in (0, 2)$, define the vector $d_t \in \qq_{\ge 0}^7$ by
		\[
		d_t \coloneqq (1, 1, 0, t, 1, 0, 0).
		\]
		The support of $d_t$ is $ABZU$, so $d_t \in M$. Moreover,
		\[
		a_i - d_t = u_i + (1 - t)e_Z
		\]
		for every $i \in \ldb 1, 3 \rdb$, where $e_Z \in \qq^7$ is the standard vector at the coordinate $Z$. The $Z$-coordinate of this difference is $2 - t > 0$, so its support equals that of $u_i$. Therefore $d_t$ divides each $a_i$. The three differences have no common base coordinate, so their only common divisor is $0$. Hence every $d_t$ is an MCD of $a_1$, $a_2$, and $a_3$. These MCDs are pairwise distinct and, since $M$ is reduced, pairwise nonassociate.
	\end{proof}
	
	The same construction also separates the D-property from the PSP property.
	
	\begin{remark}
		The tuple $(u_1, u_2, u_3)$ in the preceding proof is primitive but fails the additive PSP condition. Indeed, for each rational number $\varepsilon \in (0, 1)$, the vector $q \coloneqq -\varepsilon e_Z \in \qq^7$ does not belong to $M$, while $q + u_i \in M$ for every $i \in \ldb 1, 3 \rdb$.
	\end{remark}
	
	The monoid in Theorem~\ref{thm:monoid} cannot be realized as the reduced multiplicative monoid of a domain. Ring addition would force the existence of a vector whose support is not allowed.
	
	\begin{prop} \label{prop:nonrealizable}
		The monoid $M$ is not isomorphic to the reduced multiplicative monoid of an integral domain.
	\end{prop}
	
	\begin{proof}
		Suppose, by way of contradiction, that $M$ is isomorphic to the reduced multiplicative monoid of a domain $R$ with fraction field $K$. This isomorphism extends to the Grothendieck groups. Since $\gp(M) = \qq^7$, we may identify the group of divisibility of $R$ with $\qq^7$ and its positive cone with $M$. Let $\nu \colon K^\times \to \qq^7$ be the quotient homomorphism, and, for each $i \in \{A, B, C, Z, U, V, W\}$, let $v_i \colon K^\times \to \qq$ be its $i$th coordinate function. Membership of a nonzero element $x \in K$ in $R$ is equivalent to $\nu(x) \in M$.
		
		We first prove that $v_i$ is a valuation for every $i \in \{A, B, C, Z, U, V, W\}$. Fix such an index $i$. Take $x$ and $y$ in $K^\times$ with $x + y \ne 0$, and fix a rational number $r < \min\{v_i(x), v_i(y)\}$. Choose $h \in \qq^7$ whose $i$th coordinate is $r$ and whose other coordinates are strictly smaller than the corresponding coordinates of both $\nu(x)$ and $\nu(y)$. The vectors $\nu(x) - h$ and $\nu(y) - h$ have full positive support, so they belong to $M$. By surjectivity, there is $z \in K^\times$ with $\nu(z) = h$. Then $x/z \in R$ and $y/z \in R$, and hence $(x + y)/z \in R$. It follows that $v_i(x + y) \ge r$. Since this holds for every such $r$,
		\[
		v_i(x + y) \ge \min\{v_i(x), v_i(y)\}.
		\]
		Multiplicativity is inherited from $\nu$. The valuation inequality also gives equality whenever $v_i(x) \ne v_i(y)$.
		
		Choose $x$ and $y$ in $R \setminus \{0\}$ with
		\[
		\nu(x) = (1, 1, 0, 1, 1, 0, 0)
		\qquad \text{and} \qquad
		\nu(y) = (2, 0, 1, 2, 0, 1, 0).
		\]
		Both vectors belong to $M$. Since $\nu(x) \ne \nu(y)$ and $\nu(-1) = 0$, the sum $x + y$ is nonzero. Equality in the valuation inequality at coordinates with unequal values gives
		\[
		\nu(x + y) = (1, 0, 0, 1, 0, 0, w)
		\]
		for some $w \in \qq_{\ge 0}$. Its support is $AZ$ or $AZW$, neither of which belongs to $\mathcal{F}$. This contradicts $x + y \in R$.
	\end{proof}
	
	Proposition~\ref{prop:nonrealizable} explains why the construction does not answer the following question for domains.
	
	\begin{question} \label{ques:idf-d-mcd-finite}
		Is every IDF-domain satisfying the D-property MCD-finite?
	\end{question}
	
	By Theorem~\ref{thm:ek}, Question~\ref{ques:idf-d-mcd-finite} is equivalent to asking whether the IDF property ascends to polynomial rings over IDF-domains satisfying the D-property. Theorem~\ref{thm:tidf-d-is-psp} gives an affirmative answer under the additional Furstenberg hypothesis.

	\bigskip
	\section*{Acknowledgments}
	\label{sec:acknowledgments}
	
	This collaboration took place through PRIMES-USA, a year-long mathematics research program hosted by the MIT Department of Mathematics. The authors thank the PRIMES program for this rewarding research opportunity.

\end{document}